\documentclass[10pt]{amsart}   
\usepackage{style}
\usepackage{mathtools}
\title{Complex orientations and $\TP$ of complete DVRS}
\author{Gabriel Angelini-Knoll} 
\date{}							

\begin{document}
\begin{abstract}
Let $L$ be finite extension of $\mathbb{Q}_p$ with ring of integers $\mathcal{O}_L$. We show that periodic topological cyclic homology of $\mathcal{O}_L$, over the base $\mathbb{E}_{\infty}$-ring $\mathbb{S}_{W(\mathbb{F}_q)}[z]$ carries a $p$-height one formal group law mod $(p)$ that depends on an Eisenstein polynomial of $L$ over $\mathbb{Q}_p$ for a choice of uniformizer $\varpi\in \mathcal{O}_L$.  
\end{abstract}

\maketitle

\section{Introduction}
In Bhatt--Morrow--Scholze \cite[Prop. 11.10]{BMS19}, they compute $p$-complete periodic topological cyclic homology of the ring of integers $\mathcal{O}_L$ in a finite extension $L$ of $\mathbb{Q}_p$ over the base $\mathbb{E}_{\infty}$-ring $\bS[z]\coloneqq\Sigma^{\infty}_+\mathbb{N}$, denoted $\TP(\mathcal{O}_L/\bS[z];\mathbb{Z}_p^{\wedge})$. Specifically, they prove that $\TP_*(\mathcal{O}_L/\bS[z];\mathbb{Z}_p^{\wedge})$ is isomorphic to a periodized Frobenius twist of the Breul--Kisin module in loc. cit.\footnotemark \footnotetext{As a ring, the Frobenius twist of a Breul-Kisin module $\mathfrak{S}^{-1}$ is simply $\mathfrak{S}=W(\mathbb{F}_q)[[z]]$, but it is regarded as a $\mathfrak{S}$-algebra via the canonical Frobenius map $\varphi$, which lifts the Frobenius map on $W(\mathbb{F}_q)$ and sends $z$ to $z^p$.} The starting point of this note is the observation that the spectrum $\TP(\mathcal{O}_L/\bS[z];\mathbb{Z}_p^{\wedge})$ is complex oriented. We show that it carries a height one formal group law that is compatible with the cyclotomic structure map of $\mathbb{E}_{\infty}$-rings 
\begin{align}\label{varphi}
\varphi \colon \thinspace \TC^{-}(\mathcal{O}_L/\mathbb{S}[z];\mathbb{Z}_p^{\wedge})\to \TP(\mathcal{O}_L/\bS[z];\mathbb{Z}_p^{\wedge}),
\end{align}
from \cite[Const. 11.5]{BMS19}.
This gives evidence that periodic topological cyclic homology shifts chromatic complexity up by one, consistent with the red-shift philosophy of Ausoni--Rognes \cite{AR08}, where chromatic complexity is measured using complex orientations. Since algebraic K-theory is rarely complex oriented, it is less common to measure red-shift behavior using complex orientations, but here we make the case that it can be used as a method for studying red-shift behavior in $\TP$ over a choice of base $\mathbb{E}_{\infty}$-ring, which in this case is given by a choice of uniformizer in $\mathcal{O}_L$. 

Throughout, let $L$ denote a finite extension of $\mathbb{Q}_p$ with ring of integers $\mathcal{O}_L$, which is a complete discrete valuation ring of mixed characteristic $(0,p)$. We fix a uniformizer $\varpi_L$ in $\mathcal{O}_L$ and write $\mathbb{F}_q=\mathcal{O}_L/(\varpi_L)$ for the residue field, which is a finite field of characteristic $p$. Let $E_L(z)$ denote an Eisenstein polynomial of degree $e_L$ with leading coefficient $\mu$. This is determined (up to multiplication by a unit) by our choice of uniformizer in the sense that the quotient map $W(\mathbb{F}_q)[z]\longrightarrow \mathcal{O}_L$ sends $z$ to $\varpi_L$ and the kernel is $(E_L(z))$. This specifies an isomorphism of $W(\mathbb{F}_q)$-algebras $\mathcal{O}_L\cong W(\mathbb{F}_q)[z]/(E_L(z))$ by the classification of local fields of mixed characteristic $(0,p)$ \cite[II. \S 5 Thm. 4]{Ser79}.  
 
We write $\bS_{W(\mathbb{F}_q)}$ for the spherical Witt vectors, which is an $\mathbb{E}_{\infty}$-ring constructed in \cite[Ex. 5.2.7]{Lur21}. The choice of uniformizer $\varpi_L$ determines a map 
$\mathbb{S}_{W(\mathbb{F}_q)}[z]\to \mathcal{O}_L,$
of $\mathbb{E}_{\infty}$-rings sending $z$ to $\varpi_L$. We define
$\THH(\mathcal{O}_L/\bS_{W(\mathbb{F}_q)}[z])$  
by replacing all tensors over $\mathbb{S}$ in the standard cyclic bar construction with tensors over $\bS_{W(\mathbb{F}_q)}[z]$ (cf. Eq. \eqref{cyclic bar}).  
By \cite[Rem. 2.12]{LW20}, this is simply the $p$-completion of $\THH(\mathcal{O}_L/\mathbb{S}[z])$. We define 
\begin{align*}
\TC^{-}(\mathcal{O}_L/\mathbb{S}_{W(\mathbb{F}_q)}[z])\coloneqq&\THH(\mathcal{O}_L/\mathbb{S}_{W(\mathbb{F}_q)}[z])^{h\mathbb{T}} \text{ and } \\
\TP(\mathcal{O}_L/\mathbb{S}_{W(\mathbb{F}_q)}[z])\coloneqq&\THH(\mathcal{O}_L/\mathbb{S}_{W(\mathbb{F}_q)}[z])^{t\mathbb{T}}
\end{align*}
where $\mathbb{T}\subset \mathbb{C}^{\times}$ is the subgroup of unit vectors. As we recall in Proposition \ref{main prop}, it is shown in \cite[Prop. 11.10]{BMS19} (cf. Theorem 2.15 \cite{LW20}) that 
\[ \TC_0^{-}(\mathcal{O}_L/\bS_{W(\mathbb{F}_q)}[z]) \cong \TP_0(\mathcal{O}_L/\bS_{W(\mathbb{F}_q)}[z])\cong W(\mathbb{F}_q)[[z]]\] 
and the Tate valued Frobenius map induces a map of $W(\mathbb{F}_q)$-algebras 
\[\varphi \colon \thinspace \TC_0^{-}(\mathcal{O}_L/\bS_{W(\mathbb{F}_q)}[z]) \to \TP_0(\mathcal{O}_L/\bS_{W(\mathbb{F}_q)}[z])\]
extending the Frobenius map on $W(\mathbb{F}_q)$ by sending $z$ to $z^p$.

We say a formal group law $F$ over a commutative $\mathbb{F}_p$-algebra $R$ has \emph{$p$-height $n$} if the $p$-series satisfies 
\[ [p]_F(X)=aX^{p^n}+O(X^{p^n+1})\in R[[x]]\]  
where $a\ne 0$ in $R$.\footnote{Here, we use the definition of height from M. Lazard \cite[p.266]{Laz55} (cf. \cite[p.27]{Fro68}). This ensures that each formal group has a well-defined height. In other references, such as \cite[Rem. 4.5.7]{Smi11} this notion is called height $\ge n$. We also expect that $F$ has height (exactly) $n$ in the sense of \cite[Prop. 4.5.3]{Smi11}. In other words, we expect that $\overline{\alpha}$ is a unit for any choice of complex orientation. See Remark \ref{obstruction} for further discussion.}
Recall that $\MU_*\cong \mathbb{Z}[x_i: i\ge 1]$ and throughout let $v_n\coloneqq x_{p^n-1}$.
The spectrum $\TC^{-}(\mathcal{O}_L/\mathbb{S}_{W(\mathbb{F}_q)}[z])$ is equipped with the skeletal filtration (cf. Eq. \eqref{skeletal filtration}). We will write $x=y+O(t^2)$ if $x=y$ in the quotient $\TC_*^{-}(\mathcal{O}_L/\mathbb{S}_{W(\mathbb{F}_q)}[z])/F_4$ using the notation from Eq. \eqref{skeletal filtration}. We also write  $x= y+O(\varphi(t^2))$ in $\TP_*(\mathcal{O}_L/\mathbb{S}_{W(\mathbb{F}_q)}[z])$ if $x=y$ in the quotient of the composite
\[ 
	\begin{tikzcd}
		F_4\subset \TC_*^{-}(\mathcal{O}_L/\mathbb{S}_{W(\mathbb{F}_q)}[z])\ar[r,"\varphi"] & \TP_*(\mathcal{O}_L/\mathbb{S}_{W(\mathbb{F}_q)}[z]).
	\end{tikzcd}
\] 
See Proposition \ref{main prop} for the naming convention of classes in $\TC_*^{-}(\mathcal{O}_L/\mathbb{S}_{W(\mathbb{F}_q)}[z])$ and $\TP_*(\mathcal{O}_L/\mathbb{S}_{W(\mathbb{F}_q)}[z])$. 
\begin{thm}\label{main thm}
Let $p\ge 3$. Any choice of complex orientation 
\[ \MU\to \TC^{-}(\mathcal{O}_L/\mathbb{S}_{W(\mathbb{F}_q)}[z]))\]
induces a map on homotopy groups sending the class $v_1\in \pi_{2p-2}\MU$ to a class $\overline{v}_1=\alpha E_L(z)+O(t^2)$ for some $\alpha \in \mathcal{O}_L$ with non-trivial mod $(p)$-reduction. The spectrum $\mathbb{S}/p\otimes \TC^{-}(\mathcal{O}_L/\mathbb{S}_{W(\mathbb{F}_q)}[z])$ therefore carries a $p$-height one graded formal group law $F$ with $p$-series 
\[[p]_{F}(X)=\overline{\alpha}\overline{\mu}z^{e_L}x^{p-1}X^{p}+O(X^{p+1})+O(t^2) \]
where $\overline{\alpha}\overline{\mu}z^{e_L}\not = 0$ is the mod $(p)$ reduction of $\alpha E_L(z)$.
\end{thm}
 
\begin{cor}\label{main cor}
Let $p\ge 3$. Any complex orientation 
\[ \MU\to \TP(\mathcal{O}_L/\mathbb{S}_{W(\mathbb{F}_q)}[z]))\]
that factors through the map \eqref{varphi} induces a map on homotopy groups sending the class $v_1\in \pi_{2p-2}\MU$ to a class $\overline{v}_1^{\prime}\coloneqq \varphi(\alpha) \varphi(E_L(z)) \sigma^{p-1} + O(\varphi(t^2))$
for some $\alpha \in \mathcal{O}_L$ with non-trivial mod $(p)$-reduction $\overline{\alpha}$. The complex oriented spectrum $\mathbb{S}/p\otimes \TP(\mathcal{O}_L/\mathbb{S}_{W(\mathbb{F}_q)}[z])$ therefore carries a $p$-height one graded formal group law $G$ with $p$-series 
\[[p]_{G}(X)=\varphi(\overline{\alpha}\overline{\mu})z^{pe_L} \sigma^{p-1}X^{p}+O(X^{p+1})+O(\varphi(t^2)) \]
where $\varphi(\overline{\alpha}\overline{\mu})z^{pe_L}\not = 0$ is the mod $(p)$ reduction of $\varphi(\alpha E_L(z))$.
\end{cor}
\begin{rem2}
By \cite[Prop. 11.10]{BMS19} (cf. \cite[Thm 2.15]{LW20}), we know that $E_L(z)$ maps to $E_L(z)=0\in \pi_0\THH(\mathcal{O}_L/\mathbb{S}_{W(\mathbb{F}_q)}[z])$, so the multiplicative formal group law of Theorem \ref{main thm} reduces to the additive formal group law on $\mathbb{S}/p\otimes \THH(\mathcal{O}_L/\mathbb{S}_{W(\mathbb{F}_q)}[z])$.  Similarly, the element $E_L(z)$ maps to $E_L(\varpi_L)=0\in\TP_0(\mathbb{F}_q)$ so this formal group law also reduces to the additive formal group law on $\mathbb{S}/p\otimes \TP(\mathbb{F}_q)$. This is consistent with the red-shift philosophy.\end{rem2}
\begin{rem2}
We suspect that $\alpha$ is in fact a unit and can be chosen to be $1$ for  the right choice of complex orientation. See Remark \ref{obstruction} for further discussion. As observed in Proposition \ref{main prop}, the element $\overline{v}_1=vx^p+O(t^2)=E_L(z)x^{p-1}+O(t^2)$ maps to $\varphi(E_L(z))\sigma^{p-1}+O(\varphi(t^2))$ via both the canonical map as well as the Tate valued Frobenius map and therefore it produces a class $v_1\in \pi_{2p^2-2}\TC^{-}(\mathcal{O}_L/\mathbb{S}_{W(\mathbb{F}_q)}[z])$.
If $\alpha=1$, then the relationship between $\overline{v}_1^{\prime}$ and the element $\beta$ defined in \cite[Theorem 1.1]{LW20} is given by the formula
\[
	\overline{v}_1^{\prime} = \begin{cases}
				\beta^{p-1} & \text{ if } [L(\zeta_p), L]=1 \\ 	
				\beta & \text{ if } [L(\zeta_p), L]=p-1.\footnotemark 				\end{cases}
\]
\footnotetext{{Our results were proven independently and largely in parallel to \cite{LW20}. We have updated this paper to include references to the results of loc. cit. in their full strength as they have become available to the author.}
}
 
\end{rem2}

\subsection{Conventions}
To fix models, we work in the symmetric monoidal $\infty$-category of spectra throughout and write $\otimes$ for the symmetric monoidal product, $\mathbb{S}$ for the symmetric monoidal unit, and $\oplus$ for the coproduct. We also write $\otimes_{R}$ for the symmetric monoidal product in graded $R$-modules for a commutative ring and we write $\oplus$ for the direct sum of graded abelian groups and the distinction is clear from context. When $R=\mathbb{F}_p$, we simply write $\otimes$. We write $\mathcal{A}$ for the dual Steenrod algebra defined as $\mathcal{A}=\pi_*(\mathbb{F}_p\otimes \mathbb{F}_p)$ using our conventions. By an $\mathbb{E}_n$-ring, we mean an algebra over the little cubes $\infty$-operad $\mathbb{E}_n^{\otimes}$ in the symmetric monoidal $\infty$-category of spectra. Given a spectral sequence $E_r^{*,*}\implies \pi_*X$ with differential $d_r^{*,*}$ we say $x\in E_r^{*,*}$ is an infinite cycle if $d_s(x)=0$ for all $s\ge r$ and we say $x\in E_r^{*,*}$ is a permanent cycle if it is an infinite cycle and it is also not a boundary of a differential of any length. We write $\Lambda_{R}(x_1,\dots x_n)$ for an exterior algebra on generators $x_1,\dots x_n$ where $n$ is a positive integer or $\infty$. We write $E(n)=\Lambda_{\mathbb{F}_p}(Q_0,\dots ,Q_n)$ for the exterior algebra on the Milnor primitives $Q_0,\dots ,Q_n$  in the Steenrod algebra $\mathcal{A}^*$.

\subsection{Acknowledgements}
The author would like to thank Mona Merling, Thomas Nikolaus, and Holger Reich for feedback on an early draft of this note, and Jack Morava and Andrew Salch for inspirational conversations. The author would also like to thank an anonymous referee for their suggestions, which have led to improvements of the paper. The author was informed that similar calculations where done by Bhatt--Morrow--Scholze, Hesselholt, and Krause--Nikolaus independently in unpublished work, so the author acknowledges these other computations. 
\section{The approximate fixed point spectral sequence}
Fix a ring spectrum $\mathsf{E}$, an $\mathbb{E}_{\infty}$-ring $A$ and an $\mathbb{E}_1$-$A$-algebra $B$. The topological Hochschild homology $\THH(B/A)$ is defined as the realization of the cyclic bar construction 
\begin{align}\label{cyclic bar}
	\THH(B/A)=|\xymatrix{B \ar[r]& \ar@<1ex>[l]\ar@<-1ex>[l]B\otimes_{A}B\ar@<.5ex>[r] \ar@<-.5ex>[r]& \ar@<-1ex>[l]\ar@<1ex>[l]\ar[l] B\otimes_{A}B\otimes_{A}B \ar@<1ex>[r]\ar@<-1ex>[r] \ar[r] & \ar@<.5ex>[l]\ar@<1.5ex>[l]\ar@<-1.5ex>[l]\ar@<-.5ex>[l]\dots }|
\end{align} 
in $A$-modules and it has an action of the circle group $\mathbb{T}$. We write 
\[\TC^{-}(B/A)[k]\coloneqq F (S(\mathbb{C}^{k})_+,THH(B/A))^{\mathbb{T}}\]
and note that $\lim_k \TC^{-}(B/A)[k]=\TC^{-}(B/A)$. There is a strongly convergent, multiplicative \emph{approximate homotopy fixed point spectral sequence} 
\begin{align}\label{approx hom fix}
 	\left(E_2^{*,*}\right) [k] =\mathbb{Z}[t]/t^{k} \otimes_{\mathbb{Z}} \mathsf{E}_*(\THH(B/A))\implies \mathsf{E}_{*}(\TC^{-}(B/A)[k])
\end{align}
with length $r$ differential denoted $d_{r}^{[k]} \colon \thinspace  \left(E_2^{*,*}\right) [k] \longrightarrow \left(E_2^{*-r,*+r-1}\right) [k] $. Here $|t|=(-2,0)$ and $|x|=(0,q)$ for $x\in  \mathsf{E}_q(\THH(B/A))$ and we grade the spectral sequence using the homological Serre convention $\left (E_r^{p,q} \right )[k] \implies  \mathsf{E}_{p+q}(\TC^{-}(B/A)[k])$. 
We say that $x\in \left (E_{\infty}^{*,*}\right ) [k]$ has skeletal filtration $2m$ for some $0\le m \le k-1$ if $x$ is $t^m$-divisible. 
This provides a filtration 
\begin{align}\label{skeletal filtration}
0 \subset F_{2k-2}\subset \dots \subset F_{1} \subset  F_0
\end{align}
of $F_0=\mathsf{E}_*(\TC^{-}(B/A)[k])$ with  $F_i=F_{i+1}$ whenever $i$ is odd such that 
$\left(E_{\infty}^{*,*}\right) [k]=\bigoplus_{i=0}^{2k}F_{2i}/F_{2i+2}$
where $F_{2i}/F_{2i+2}=\pi_*(\Sigma^{2i}THH(B/A))$. 
We write 
\begin{align}\label{theta}
	\theta \colon \thinspace F_{2m} \longrightarrow F_{2m}/F_{2m+2}
\end{align}
for the canonical quotient map were the non-negative integer $m$ is understood from context. Therefore, given $x\in F_{2m}\subset \mathsf{E}_*(\TC^{-}(B/A)[k])$, there is a corresponding element $\theta(x)\in F_{2m}/F_{2m+2}\subset \left (E_{\infty}^{*,*}\right )[k]$. 

We compute $\pi_*\TC^{-}(\mathcal{O}_L/\mathbb{S}_{W(\mathbb{F}_q)}[z])[2]$,
which follows directly from \cite{BMS19,KN19,LW20}, but we provide a proof for completeness. In \cite[Thm. 3.1]{KN19}, they compute 
\[ \THH_*(\mathcal{O}_L/\mathbb{S}_{W(\mathbb{F}_q)}[z])\cong \mathcal{O}_L[x]\]
where $\mathcal{O}_L[x]$ is a graded polynomial algebra over $\mathcal{O}_L$ with $|x|=2$. At the outset, we fix a preferred choice of indecomposable algebra generator $x\in \THH_2(\mathcal{O}_L/\mathbb{S}_{W(\mathbb{F}_q)}[z])$ by defining $x\coloneqq\sigma^{2}E_L(z)$ where 
\[ \sigma^{2}\colon \thinspace \Sigma \left ( \mathcal{O}_{L}/\mathbb{S}_{W(\mathbb{F}_q)}[z]  \right )\to \THH(\mathcal{O}_L/\bS_{W(\mathbb{F}_q)}[z])\]
is the map constructed in \cite[Construction A.1.2, Example A.2.4]{HW20} and $\mathcal{O}_{L}/\mathbb{S}_{W(\mathbb{F}_q)}[z]$ denotes the cofiber of the unit map $\mathbb{S}_{W(\mathbb{F}_q)}[z]\to \mathcal{O}_{L}$. Note that this choice depends on a choice of Eisenstein polynomial. Let $\mathfrak{S}\coloneqq W(\mathbb{F}_q)[[z]]$ throughout this section.
\begin{lem}\label{easy lem}
There is an isomorphism of graded rings 
\[ \pi_*(\TC^{-}(\mathcal{O}_L/\mathbb{S}_{W(\mathbb{F}_q)}[z])[2])\cong \mathfrak{S}[x,t]/(tx-E_L(z),t^2)\]	
with $|z|=0$, $|x|=2$, and $|t|=-2$, 
\end{lem}
\begin{proof}
The spectral sequence \eqref{approx hom fix} for $k=2$ collapses at the $E_2$-page for bi-degree reasons. Note that $x=\theta(\widetilde{x})$ for some choice of lift $\widetilde{x}\in \pi_{2}(\TC^{-}(\mathcal{O}_L/\mathbb{S}_{W(\mathbb{F}_q)}[z])$ along the surjection
\begin{align}\label{specific theta}
	\theta \colon \thinspace \pi_{2}(TC^{-}(\mathcal{O}_L/\mathbb{S}_{W(\mathbb{F}_q)}[z])\longrightarrow \pi_2(THH(\mathcal{O}_L/\mathbb{S}_{W(\mathbb{F}_q)}[z])
\end{align}
defined in \eqref{theta}. We choose a lift $\widetilde{x}$ to be the generator of $\pi_{2}(TC^{-}(\mathcal{O}_L/\mathbb{S}_{W(\mathbb{F}_q)}[z])$ in the complement $F_0\backslash F_1$ as an $\mathfrak{S}$-module and call this choice  of lift simply $x$ by abuse of notation. Note that at this point of the argument, we have not yet resolved the extensions, so a priori we only know that there is an exact sequence
\[ 
	\begin{tikzcd}
		0 \ar[r] & \mathcal{O}_L\{tx^2\} \ar[r] &  \pi_{2}(TC^{-}(\mathcal{O}_L/\mathbb{S}_{W(\mathbb{F}_q)}[z])\ar[r,"\theta"] &  \mathcal{O}_L\{x\}\ar[r] &0
	\end{tikzcd}
\]
Nonetheless, this does not affect our choice of lift $\widetilde{x}$.\footnote{This choice of lift is compatible with a choice of lift of the canonical quotient 
$\theta \colon \thinspace W(\mathbb{F}_q)[[z]]\{x\}\longrightarrow \mathcal{O}_L\{x\}$ in the the homotopy fixed point spectral sequence, where now this corresponds to the canonical quotient  $\theta \colon \thinspace \TC_2^{-}(\mathcal{O}_L/\mathbb{S}_{W(\mathbb{F}_q)}[z])\longrightarrow THH_2(\mathcal{O}_L) \cong W(\mathbb{F}_q)[[z]]/(E_L(z))$.}

We compute a hidden multiplicative extension $E_L(z)=t\sigma^2E_L(z)=tx$ by applying \cite[Lemma A.4.1]{HW20}.\footnote{Note that $\TC^{-}(\mathcal{O}_L/\mathbb{S}_{W(\mathbb{F}_q)}[z])[2]=\lim_{\mathbb{C}P^1}\HH(\mathcal{O}_L/\mathbb{S}_{W(\mathbb{F}_q)}[z])$ in the notation from \cite[Lemma A.4.1]{HW20}} Note that there is no indeterminacy in the choice of $t$ here because it is in highest skeletal filtration in the approximate homotopy fixed point spectral sequence
We claim that there are is no further room for multiplicative extensions besides those produced by the identification $E_L(z)=tx$. To see this, note that the remaining products that are zero in the $E_{\infty}$-page are in skeletal filtration two or higher. Since multiplicative extensions must raise skeletal filtration in the approximate fixed point spectral sequence and the spectral sequence is concentrated in skeletal filtration zero and two, there is no room for further hidden multiplicative extensions. 
\end{proof}
\noindent Note that $\THH(\mathcal{O}_L/\mathbb{S}_{W(\mathbb{F}_q)}[z])$ is an $\mathbb{E}_\infty$-ring cyclotomic spectrum by \cite[Const. 11.5]{BMS19}. We also record the more general computation of $\TC_*^{-}(\mathcal{O}_L/\mathbb{S}_{W(\mathbb{F}_q)}[z])$ and $\TP_*(\mathcal{O}_L/\mathbb{S}_{W(\mathbb{F}_q)}[z])$. 
\begin{proposition}[Proposition 11.10 \cite{BMS19} (cf. Theorem 2.15 \cite{LW20})]\label{main prop}
There are isomorphisms 
\begin{align*}
\TC_*^{-}(\mathcal{O}_L/\mathbb{S}_{W(\mathbb{F}_q)}[z])\cong & \mathfrak{S}[x,v]/(vx-E_L(z))\text{ and }\\
\TP_*(\mathcal{O}_L/\mathbb{S}_{W(\mathbb{F}_q)}[z])\cong & \mathfrak{S}[\sigma,\sigma^{-1}] 
\end{align*}
of graded rings with $|x|=2$, $|v|=-2$, and $|z|=0$. The canonical map 
\[\operatorname{can} \colon \thinspace \TC^{-}_*(\mathcal{O}_L/\mathbb{S}_{W(\mathbb{F}_q)}[z])\to \TP_*(\mathcal{O}_L/\mathbb{S}_{W(\mathbb{F}_q)}[z])\]
sends u to $E_L(z)\sigma$ and $v$ to $\sigma^{-1}$. 

The Frobenius map 
\[\varphi \colon \thinspace \TC^{-}_0(\mathcal{O}_L/\mathbb{S}_{W(\mathbb{F}_q)}[z])\to \TP_0(\mathcal{O}_L/\mathbb{S}_{W(\mathbb{F}_q)}[z])\]
is the lift of Frobenius on $W(\mathbb{F}_q)$ and sends $z$ to $z^p$. More generally, the map 
\begin{align}\label{frob} 
\varphi \colon \thinspace \TC^{-}(\mathcal{O}_L/\mathbb{S}_{W(\mathbb{F}_q)}[z])\to \TP(\mathcal{O}_L/\mathbb{S}_{W(\mathbb{F}_q)}[z])
\end{align}
induces a map of graded rings on homotopy groups $\pi_*$ by sending $x$ to $\sigma$, $v$ to $\varphi(E_L(z))\sigma^{-1}$ 
\end{proposition}
\begin{rem2}
Note that we have chosen to use the notation from \cite[Proposition 11.10]{BMS19} except for the element $x$. The element $v$ may therefore differ from the element $t$ in the $E_2$-page of the homotopy fixed point spectral sequence, but we still know that $v=t+O(t^2)$. Similarly, the element $\sigma^{-1}$ may differ from the element $t$ in the $E_2$-page of the Tate spectral sequence, but we still know that $\sigma^{-1}=t+O(t^2)$. 
\end{rem2}

\section{Complex orientations and periodic topological cyclic homology} 
This section is devoted to proving Theorem \ref{main thm}.
\begin{proposition}\label{homological homotopy fixed point ss}
Let $\left (E_{r}^{*,*}\right )[2]$ be the $E_{r}$-page of the spectral sequence \eqref{approx hom fix} for $\mathsf{E}=\mathbb{F}_p$, $A=\mathbb{S}_{W(\mathbb{F}_q)}[z]$ and $B=\mathcal{O}_L$. 
\begin{comment}
Then we identify
$\left (E_{\infty}^{*,*}\right )[2]= M \otimes \mathbb{F}_q[\overline{\xi}_i : i\ge 1]   $
as $\cA$-comodule algebras where 
\[ M\cong \Lambda_{\mathbb{F}_p}(\tau_k^{\prime} : k\ge 2) \otimes \mathbb{F}_p[x]/x^p\otimes \mathbb{F}_p[t]/t^2\oplus \mathbb{F}_p\left \{ t x^{\ell-1} \overline{\tau}_1 ,x^k : k,\ell \ge 1 \right \}\]
and $\tau_k^{\prime}=\overline{\tau}_k-x^{p^k-{p^{k-1}}}\overline{\tau}_{k-1}$. 
In particular, we conclude that 
\end{comment}
There is an isomorphism of $\mathcal{O}_L/(p)$-modules 
\[
	\pi_{s}\left ( \mathbb{F}_p\otimes \TC^{-}(\mathcal{O}_L/\mathbb{S}_{W(\mathbb{F}_q)}[z])[2] \right ) \cong 
		\begin{cases}
			\mathcal{O}_L/(p)\{\theta(t) \}  & \text{ if } s=-2, \\
			\mathcal{O}_L/(p)\{\theta(x)^j,\theta(t) \theta(x)^{j+1} \} & \text{ if } s=2j \text{ and }0\le  2j<2p-4, \\
			\mathcal{O}_L/(p)\{\theta(t) \theta(x)^{p-1},\theta(t) \theta( \overline{\xi}_1) \} & \text{ if } s=2p-4, \\
			\mathcal{O}_L/(p)\{\theta (t\overline{\tau}_1)\} & \text{ if } s=2p-3, \\
			\mathcal{O}_L/(p)\{\theta(x)^{p-1},\theta(\overline{\xi}_1)\} & \text{ if } s=2p-2, \\
			\mathcal{O}_L/(p)\{\theta(x) (\theta( t\overline{\tau}_1) \} & \text{ if } s=2p-1, \\
			\mathcal{O}_L/(p)\{\theta(x)^p\} & \text{ if } s=2p \\
			0 & \text{ otherwise if } s\le 2p 
		\end{cases}
\]
where $\theta$ is defined as in \eqref{theta}. 
The $\mathcal{A}$-co-action on  $\theta(\overline{\xi}_1),\theta(x)^{p-1}\in \pi_{2p-2}\left ( \mathbb{F}_p\otimes \TC^{-}(\mathcal{O}_L/\mathbb{S}_{W(\mathbb{F}_q)}[z])[2] \right )$ and $\theta (x ) \theta ( t\overline{\tau}_1) \in \pi_{2p-1}\left ( \mathbb{F}_p\otimes \TC^{-}(\mathcal{O}_L/\mathbb{S}_{W(\mathbb{F}_q)}[z])[2] \right )$
\begin{align*}
	\psi(\theta (\overline{\xi}_1) )=&\theta (\xi_1)\otimes 1+1\otimes \theta ( \overline{\xi}_1), \\
	\psi( \theta(x)^{p-1})=&1\otimes \theta(x)^{p-1}, \text{ and }\\
	\psi(\theta (x) \theta( t\overline{\tau}_1))=&1\otimes \theta (x ) \theta ( t \overline{\tau}_1)+ \overline{\tau}_1\otimes \theta (t ) \theta (x) + \overline{\tau}_0\otimes \theta (t) \theta  (\overline{\xi}_1).
\end{align*}
respectively. 
\end{proposition}

\begin{proof} 
The spectral sequence \eqref{approx hom fix} for $E=\mathbb{F}_p$, $A=\mathbb{S}_{W(\mathbb{F}_q)}[z]$ and $B=\mathcal{O}_L$  and $k=2$ has signature
\[ \left (E_2^{*,*}\right )[2]=\mathbb{Z}[t]/t^2\otimes_{\mathbb{Z}} (\mathbb{F}_p)_*(\THH(\mathcal{O}_L/\mathbb{S}_{W(\mathbb{F}_q)}[z])))\implies  \pi_{*}\left ( \mathbb{F}_p\otimes \TC^{-}(\mathcal{O}_L/\mathbb{S}_{W(\mathbb{F}_q)}[z])[2] \right ). \]
Since $\THH(\mathcal{O}_L/\mathbb{S}_{W(\mathbb{F}_q)}[z]))$ is an $E_{\infty}$-$\mathcal{O}_L$-algebra which is free as an $\mathcal{O}_L$-module, we compute 
\[ \mathbb{Z}[t]/t^2\otimes_{\mathbb{Z}}  \pi_* ( \mathbb{F}_p\otimes\THH(\mathcal{O}_L/\mathbb{S}_{W(\mathbb{F}_q)}[z])) =  \mathbb{F}_p[t]/t^2 \otimes \mathcal{O}_L/(p) \otimes \mathcal{A}//E(0)_* \otimes \mathbb{F}_p[x]\]
regarded as a $\mathcal{A}$-comodule algebra with the evident algebra structure, the $\mathcal{A}$-co-action on $(\mathbb{F}_p)_*(\mathcal{O}_L)$ is given by restriction along the inclusion 
$(\mathbb{F}_p)_*(\mathcal{O}_L)\rightarrow (\mathbb{F}_p)_*(\mathcal{O}_L/(p))$ where $(\mathbb{F}_p)_*( \mathcal{O}_L/(p))\cong \mathcal{A}\otimes  \mathcal{O}_L/(p)$ has co-action induced by the co-product, $x$ is a $\mathcal{A}$-comodule primitive, and $\mathbb{F}_p[t]/t^2$ has the $\mathcal{A}$-co-action of $H^*(\mathbb{C}P^2;\mathbb{F}_p)$ computed in \cite{Mil58}. 	
We claim the first nontrivial differentials are 
\[ d_{2}^{[2]}(\bar{\tau}_k)=tx^{p^k}\]
for $k\ge 1$.  By \cite[Prop. 3.2]{BR05} it suffices to show that $\sigma \bar{\tau}_k=x^{p^k}$. 
We then observe that the map 
\[\THH(\mathcal{O}_L/\bS_{W(\mathbb{F}_q)}[z])\to \THH(\mathbb{F}_q)\] 
is a map of $\mathbb{E}_{\infty}$-rings and therefore the map in mod $p$ homology is compatible with Dyer-Lashof operations. 
Consequently, we can use the formula $Q^{p^k}\sigma=\sigma Q^{p^k}$ by \cite{Bok87a} and the computation
\[\sigma \overline{\tau}_k=\sigma Q^{p^k}\overline{\tau}_0=Q^{p^k}(\sigma \tau_0)=(\mu_0)^{p^k} \]
by \cite{BMMS86} for $k\ge 1$, to show that $\sigma(\overline{\tau}_k)=x^{p^k}$ in $(\mathbb{F}_p)_*(\THH(\mathcal{O}_L/\bS_{W(\mathbb{F}_q)}[z]))$. The spectral sequence collapses at the $E_4$-page for bi-degree reasons.
We claim that there are no hidden $\mathcal{A}$-comodule extensions involving the classes 
\begin{align*}
\theta (\overline{\xi}_1),\theta( x)^{p-1}\in & \pi_{2p-2}\left ( \mathbb{F}_p\otimes \TC^{-}(\mathcal{O}_L/\mathbb{S}_{W(\mathbb{F}_q)}[z])[2] \right )\text{ and }\\ 
\theta (x) \theta (t \overline{\tau}_1) \in & \pi_{2p-1}\left ( \mathbb{F}_p\otimes \TC^{-}(\mathcal{O}_L/\mathbb{S}_{W(\mathbb{F}_q)}[z])[2] \right ).
\end{align*} 
Since the stated co-actions are the $\mathcal{A}$-co-actions on the corresponding classes  $\overline{\xi}_1,x^{p-1}\in \left( E_\infty^{*,*} \right )[2]$ and $x \cdot (t\overline{\tau}_1)\in \left( E_\infty^{*,*} \right )[2]$, this implies the stated $\mathcal{A}$-comodule actions on the abutments. The fact that $\theta (x)\theta (t\overline{\tau}_1)$ has the same $\mathcal{A}$-comodule action as $x\cdot (t\overline{\tau}_1)$ follows because there are no non-trivial elements in strictly higher skeletal filtration than $2$, which is the skeletal filtration of $\theta (x) \theta( t \overline{\tau}_1)$. 

We claim that $\theta (x)^{p-1}\in \pi_{2p-2} \left ( \mathbb{F}_p\otimes \TC^{-}(\mathcal{O}_L/\mathbb{S}_{W(\mathbb{F}_q)}[z])[2] \right )$ is a $\mathcal{A}$-comodule primitive. This follows because $\pi_{*} \left ( \mathbb{F}_p\otimes \TC^{-}(\mathcal{O}_L/\mathbb{S}_{W(\mathbb{F}_q)}[z])[2] \right )$ is an $\mathcal{A}$-comodule algebra and the class 
\[\theta(x) \in \pi_{2} \left ( \mathbb{F}_p\otimes \TC^{-}(\mathcal{O}_L/\mathbb{S}_{W(\mathbb{F}_q)}[z])[2] \right )\] 
is a $\mathcal{A}$-comodule primitive for degree reasons. Consequently, the class $\theta(x)^{p-1}$ is also an $\mathcal{A}$-comodule primitive. We then resolve the hidden $\mathcal{A}$-comodule structure on $\overline{\xi}_1$ by arguing by contradiction. For degree reasons, we know that $\mathcal{A}$-comodule structure on $\overline{\xi}_1$ is 
\begin{align*} 
	\psi(\overline{\xi}_1)=&\overline{\xi}_1\otimes 1+1\otimes \overline{\xi}_1 + W\cdot \overline{\tau}_0\otimes t\overline{\tau}_1 
\end{align*}
for some $W\in \mathcal{O}_L/p$. Suppose $W$ is non-zero. Then we consider the map of Adams spectral sequences with map of $E_1$-terms 
\[ \overline{\mathcal{A}}^{\otimes \bullet} \longrightarrow \overline{\mathcal{A}}^{\otimes \bullet} \otimes \pi_{*}\left ( \mathbb{F}_p\otimes \TC^{-}(\mathcal{O}_L/\mathbb{S}_{W(\mathbb{F}_q)}[z])[2] \right )\]
given by the map of normalized cobar complexes where $\overline{\mathcal{A}}$ is the kernel of the augmentation $\mathcal{A}\to \mathbb{F}_p$. 
In the source, the class $[\overline{\xi}_1] \in \overline{\mathcal{A}} $ is a permanent cycle, which projects onto the class $\alpha_1\in \pi_{2p-3}\mathbb{S}_p^{\wedge}$ in the abutment. It maps to a class 
\[ [\overline{\xi}_1\otimes 1] \in \overline{\mathcal{A}} \otimes \pi_{*}\left ( \mathbb{F}_p\otimes \TC^{-}(\mathcal{O}_L/\mathbb{S}_{W(\mathbb{F}_q)}[z])[2] \right )\]
and by the map of spectral sequences this element is an infinite cycle. If it is not a boundary, then it is a permanent cycle and there there is some non-trivial element in $\pi_{2p-3}  ( \TC^{-}(\mathcal{O}_L/\mathbb{S}_{W(\mathbb{F}_q)}[z])[2])$. This contradicts the computation from Lemma \ref{easy lem}. Therefore, the class $\overline{\xi}_1\otimes 1$ must be a boundary of a $d_1$ in the Adams spectral sequence with signature
\[E_1^{*,*}= \overline{\mathcal{A}}^{\otimes \bullet} \otimes \pi_{*}\left ( \mathbb{F}_p\otimes \TC^{-}(\mathcal{O}_L/\mathbb{S}_{W(\mathbb{F}_q)}[z])[2] \right ) \implies  \pi_* (\TC^{-}(\mathcal{O}_L/\mathbb{S}_{W(\mathbb{F}_q)}[z])[2] ),\]
which strongly converges to the desired abutment because $\TC^{-}(\mathcal{O}_L/\mathbb{S}_{W(\mathbb{F}_q)}[z])[2] )$ is bounded below, has finite type homology, and has $p$-complete homotopy groups. The only possible elements that could hit  $\overline{\xi}_1\otimes 1$ are $\theta (\xi_1)$ and $\theta (x)^{p-1}$ for bi-degree reasons. We already proved that $\theta (x)^{p-1}$ is an $\mathcal{A}$-comodule primitive, so $d_1(\theta (x ) ^{p-1})=0$. We can also compute
\[d_1(\widetilde{\xi}_1)=\overline{\xi}_1\otimes 1+W\cdot \overline{\tau}_0\otimes t\overline{\tau}_1.\] 
Since $W\cdot \overline{\tau}_0\otimes t\overline{\tau}_1$ is not a boundary of a $d_1$ when $W\ne 0$, we would not have a differential $d_1(\overline{\xi}_1)=\overline{\xi}_1\otimes 1$ as needed if $W\ne 0$. We conclude that $W=0$. 
\end{proof}

\begin{lem}\label{key lem}
Any choice of complex orientation 
\[ \MU\longrightarrow \TC^{-}(\mathcal{O}_L/\mathbb{S}_{W(\mathbb{F}_q)}[z])[2]\]
induces a map on homotopy sending 
$v_1$ to $\alpha \cdot E_L(z)x^{p-1}$ for some $\alpha \in \mathcal{O}_L$ that is not divisible by $p$.
\end{lem}
\begin{proof}
Note that $\TC^{-}(\mathcal{O}_L/\bS_{W(\mathbb{F}_q)}[z])[2]$ is an $\mathbb{E}_{\infty}$-ring and it is concentrated in even degrees. 
Consequently, there exists a complex orientation 
$\MU \to \TC^{-}(\mathcal{O}_L/\bS_{W(\mathbb{F}_q)}[z])[2],$
and it can be lifted to a map of $\mathbb{E}_2$-rings by \cite[Theorem 1.2.]{CM15}. For any such complex orientation, there is a multiplicative map of Adams spectral sequences with $E_1$-pages 
\begin{align}
\overline{\mathcal{A}}^{\otimes \bullet} \otimes (\mathbb{F}_p)_*(\MU))\to  \overline{\mathcal{A}}^{\otimes \bullet} \otimes (\mathbb{F}_p)_*(\TC^{-}(\mathcal{O}_L/\mathbb{S}_{W(\mathbb{F}_q)}[z])[2])).
\end{align}
The source strongly converges to $\pi_*\MU_p^{\wedge}$ because $(\mathbb{F}_p)_*(\MU)$ is finite type and $\MU$ is bounded below and the target strongly converges to $\pi_*(\TC^{-}(\mathcal{O}_L/\mathbb{S}_{W(\mathbb{F}_q)}[z])[2])$ because $\TC^{-}(\mathcal{O}_L/\mathbb{S}_{W(\mathbb{F}_q)}[z])[2]$ has finite type homology and it is $p$-complete and bounded below. 
The element 
\[a_1\coloneqq [\overline{\tau}_1\otimes 1-\overline{\tau}_0\otimes \overline{\xi}_1 ] \in \overline{\cA}_*\otimes (\mathbb{F}_p)_*(\MU) \]
is known to be a permanent cycle by \cite[3.1.10]{Rav86} with $v_1 \in \pi_{2p-2}(\MU)$ projecting onto $a_1\in E_{2p-1,1}^{\infty}$. 
The map of normalized cobar complexes 
\[\overline{\cA}_*^{\otimes \bullet}\otimes (\mathbb{F}_p)_*(\MU)\longrightarrow \overline{\cA}_*^{\otimes \bullet}\otimes (\mathbb{F}_p)_*(\TC^{-}(\mathcal{O}_L/\mathbb{S}_{W(\mathbb{F}_q)[z]})[2])\]
sends $a_1$ to 
\[a_1^{\prime}= [\overline{\tau}_1\otimes 1-\overline{\tau}_0\otimes \overline{\xi}_1 ]\in \overline{\cA}_*\otimes  (\mathbb{F}_p)_*(\TC^{-}(\mathcal{O}_L/\mathbb{S}_{W(\mathbb{F}_q)[z]})[2]).\]
This already implies that $a_1^{\prime}$ is an infinite cycle.  It suffices to prove that $a_1^{\prime}$ is not the boundary of a $d_1$-differential in the Adams spectral sequence. However, we computed that
\[ (\mathbb{F}_p)_{2p-1}(\TC^{-}(\mathcal{O}_L/\mathbb{S}_{W(\mathbb{F}_q)}[z]))=\mathbb{F}_p\{\theta (x) \theta ( t\overline{\tau}_1)\}\] 
in Lemma \ref{key lem}. Also, by Lemma \ref{key lem}, there is a differential $d_1(\theta (x) \theta ( t\overline{\tau}_1))=\overline{\tau}_1\otimes \theta (t ) \theta (x) + \overline{\tau}_0\otimes \theta (t) \theta  (\overline{\xi}_1)$, but clearly $d_1(\theta (x) \theta ( t\overline{\tau}_1))\ne a_1^{\prime}$, so since there are no other possible sources for a $d_1$-differential in this bi-degree, we conclude that $a_1$ is a permanent cycle in the Adams spectral sequence, which projects onto some non-zero class in $\pi_{2p-2} (\TC^{-}(\mathcal{O}_L/\mathbb{S}_{W(\mathbb{F}_q)[z]})[2])$. We therefore conclude that $v_1\in \pi_{2p-2}(\MU)$ maps to some $0\ne \widetilde{v}_1\in \TC^{-}_{2p-2}(\mathcal{O}_L/\mathbb{S}_{W(\mathbb{F}_q)}[z])[2]$. 
This implies that $v_1\in\pi_{2p-2}\MU$ maps non-trivially to a class $\overline{v}_1\in \pi_{2p-2}\TC^{-}(\mathcal{O}_L/\mathbb{S}_{W(\mathbb{F}_q)}[z])\cong W(\mathbb{F}_q)[[z]]$. Since $\THH(\mathcal{O}_L/\mathbb{S}_{W(\mathbb{F}_q)}[z])$ is an $\mathcal{O}_L$-module and it is complex oriented, we know that $v_1$ maps to zero in $\THH_{2p-2}(\mathcal{O}_L/\mathbb{S}_{W(\mathbb{F}_q)}[z])$. Consequently, writing 
\[ 
	\theta \colon \thinspace 
	\TC_{2p-2}^{-}(\mathcal{O}_L/\mathbb{S}_{W(\mathbb{F}_q)}[z])\cong W(\mathbb{F}_q)[[z]]\{x^{p-1}\}
	\longrightarrow \THH_{2p-2}(\mathcal{O}_L/\mathbb{S}_{W(\mathbb{F}_q)}[z])\cong \mathcal{O}_L\{x^{p-1}\}
\]
for the canonical quotient from \eqref{theta}, then $\theta(\overline{v}_1)=0 $ and $\overline{v}_1$ is detected by a class in filtration $2$. In $\TC^{-}_{2p-2}(\mathcal{O}_L/\mathbb{S}_{W(\mathbb{F}_q)}[z])[2]$, the only non-trivial element in skeletal filtration $2$ is $\alpha \cdot tx^p$ for some $0\ne \alpha \in \mathcal{O}_L$, by Lemma \ref{easy lem}. The fact that $\alpha \cdot tx^p=\alpha \cdot E_L(z)x^{p-1}$ also follows from Lemma \ref{easy lem}.

To see that  $\tilde{v}_1$ is not $p$-divisible, we apply exactly the same argument with $\MU$ replaced with $\mathbb{S}/p\otimes \MU$ and $\TC^{-}_{2p-2}(\mathcal{O}_L/\mathbb{S}_{W(\mathbb{F}_q)}[z])[2]$ replaced with $\mathbb{S}/p\otimes \TC^{-}_{2p-2}(\mathcal{O}_L/\mathbb{S}_{W(\mathbb{F}_q)}[z])[2]$. Again, we know that $a_1\in \overline{\cA}\otimes (\mathbb{F}_p)_*(\mathbb{S}/p\otimes \MU)$ is a permanent cycle in the Adams spectral sequence for $\mathbb{S}/p\otimes \MU$ and it maps to $a_1^{\prime}$. Again it suffices to check that it is not the boundary of a $d_1$ and the same argument applies for the class $\theta (x)\theta(t\overline{\tau}_1)$. Since 
\[\pi_{2p-2} (\mathbb{F}_p\otimes \mathbb{S}/p\otimes \TC^{-}_{2p-2}(\mathcal{O}_L/\mathbb{S}_{W(\mathbb{F}_q)}[z])[2])=\mathbb{F}_p\{\theta (x)\theta(t\overline{\tau}_1), \tau_0\cdot \theta (x)^{p-1},\tau_0\cdot \theta (\overline{\xi}_1\}\]
it will suffice to compute the differentials on $\tau_0\cdot \theta (x)^{p-1},\tau_0\cdot \theta (\overline{\xi}_1)$.
The computation of these differentials follow from the co-action $\psi(\tau_0)=\tau_0\otimes 1+1\otimes \tau_0$ and the computation of the $\mathcal{A}$-co-action on $\theta (x)^{p-1}$ and $\theta (\overline{\xi}_1)$ in Proposition \ref{homological homotopy fixed point ss}. Therefore, we compute 
\begin{align*} 
	d_1(\tau_0\cdot \theta (x)^{p-1})=& \tau_0 \otimes \theta (x)^{p-1}\\
	d_1(\tau_0\cdot \theta (\overline{\xi}_1))=& \tau_0 \otimes \theta (\overline{\xi}_1)+\tau_0\overline{\xi}_1\otimes 1
\end{align*}
and note that 
\[ c_1 \cdot d_1(\theta (x) \theta ( t\overline{\tau}_1))+c_2\cdot d_1(\tau_0\cdot \theta (x)^{p-1}) +c_3\cdot d_1(\tau_0\cdot \theta (\overline{\xi}_1)) \ne a_1^{\prime} .\]
for any $c_1,c_2,c_3\in \mathcal{O}_L/(p)$. Finally, observe that the same computation as that of Lemma \ref{easy lem}, but with $\mathsf{E}=\mathbb{S}$ replaced with $\mathsf{E}=\mathbb{S}/p$ implies that the only $t$-divisible class in 
$(\mathbb{S}/p)_{2p-2}(\TC^{-}(\mathcal{O}_L/\mathbb{S}_{W(\mathbb{F}_q)}[z])[2]))$ is $\overline{\alpha}\cdot tx^p$ for some $\overline{\alpha}\in \mathcal{O}_L/(p)$. Since we know there there is some non-trivial element 
\[\tilde{v}_1\in (\mathbb{S}/p)_{2p-2}(\TC^{-}(\mathcal{O}_L/\mathbb{S}_{W(\mathbb{F}_q)}[z])[2]))\] 
in skeletal filtration $2$, by the same argument as before, we know that $\overline{\alpha}\ne 0$.
\end{proof}  
\begin{proof}[Proof of Theorem \ref{main thm}]
The fact that the map $\MU_*\longrightarrow \TC_*^{-}(\mathcal{O}_L/\mathbb{S}_{W(\mathbb{F}_q)}[z])$ sends $v_1$ to 
\[\overline{v}_1=\alpha\cdot E_{L}(z)x^{p-1}+O(t^2)=\alpha vx^{p-1}+O(t^2),\] 
follows directly from Lemma \ref{key lem} since the $\mathbb{E}_{\infty}$-ring map 
$\MU_*\longrightarrow \TC_*^{-}(\mathcal{O}_L/\mathbb{S}_{W(\mathbb{F}_q)}[z])[2]$ factors through the $\mathbb{E}_{\infty}$-ring map $\MU_*\longrightarrow \TC_*^{-}(\mathcal{O}_L/\mathbb{S}_{W(\mathbb{F}_q)}[z])$ and the map  $\TC_*^{-}(\mathcal{O}_L/\mathbb{S}_{W(\mathbb{F}_q)}[z])\longrightarrow \TC_*^{-}(\mathcal{O}_L/\mathbb{S}_{W(\mathbb{F}_q)}[z])[2]$ sends classes of the same name to classes of the same name modulo skeletal filtration strictly greater than $2$. 
Applying Proposition \ref{main prop}, we know that $v_1$ maps to $\varphi(\alpha\cdot E_{L}(z))\sigma^{p-1})$ in $\TP_*(\mathcal{O}_L/\mathbb{S}_{W(\mathbb{F}_q)}[z])$ modulo the Frobenius image of elements in skeletal filtration $2$ or higher. Since the composite map 
\[
	\begin{tikzcd}
		\MU_*\ar[r] & \TC^{-}_*(\mathcal{O}_L/\mathbb{S}_{W(\mathbb{F}_q)}[z]))\ar[r,"\varphi"] & \TP_*(\mathcal{O}_L/\mathbb{S}_{W(\mathbb{F}_q)}[z])
	\end{tikzcd}
\] 
is a ring map and $\TP_*(\mathcal{O}_L/\mathbb{S}_{W(\mathbb{F}_q)}[z])$ is polynomial by Proposition \ref{main prop}, we know that the element $\varphi(\alpha)\cdot \varphi(E_{L}(z))\sigma^{p-1}$ is $v_1$-periodic. 
Consequently, we know that $\overline{v}_1$ is $v_1$-periodic in $\TC^{-}_*(\mathcal{O}_L/\mathbb{S}_{W(\mathbb{F}_q)}[z])$ as well. We can therefore compute that the graded formal group law $F$ of $\mathbb{S}/p\otimes \TC^{-}(\mathcal{O}_L/\mathbb{S}_{W(\mathbb{F}_q)}[z])$ has $p$-series 
\[ [p]_F(X)= \overline{\alpha}\overline{\mu}\cdot z^{e_L}X^{p}+O(X^{p+1})+O(t^2)\]
where $\overline{\alpha}\overline{\mu}z^{e_L}\ne 0$ is the mod $(p)$-reduction of $\alpha \cdot E_L(z)$. 
\end{proof}
\begin{proof}[Proof of Corollary \ref{main cor}]
As in the proof of Theorem \ref{main thm}, we know that the image of $v_1$ is $\overline{v}_1^{\prime}=\varphi(\alpha E_L(z))\sigma^{p-1}$ modulo the Frobenius image of elements in skeletal filtration greater or equal to $4$. We also already observed that $\overline{v}_1^{\prime}$ is $v_1$-periodic because $\TP_*(\mathcal{O}_L/\mathbb{S}_{W(\mathbb{F}_q)}[z])$ is polynomial. Passing to the mod $(p)$ reduction, we therefore have 
\[ [p]_G(X)= \overline{\alpha}\overline{\mu}\cdot z^{e_L}X^{p}+O(X^{p+1})+O(\varphi(t^2))\]
where $G$ is the $p$-typical formal group law of $\mathbb{S}/p\otimes \TP(\mathcal{O}_L/\mathbb{S}_{W(\mathbb{F}_q)}[z])$. 
\end{proof}

\begin{rem2}\label{obstruction}
We expect that $\alpha$ is always a unit and that it can be chosen to be $1$ using the right choice of complex orientation. One could try prove that $\alpha$ is a unit, by replacing $\mathbb{S}/p$ with the Moore spectrum $\mathbb{S}\mathbb{F}_q$ of $\mathbb{F}_q$ in the proof of Lemma \ref{key lem}. To do this, one would need a ring spectrum structure on the smash product $\mathbb{S}\mathbb{F}_q\otimes R$ of the Moore spectrum $\mathbb{S}\mathbb{F}_q$ with $R$ for $R\in \{\MU,\TC^{-}(\mathcal{O}_L/\mathbb{S}_{W(\mathbb{F}_q)}[z])[2]\}$ so that the $d_1$-differentials in the Adams spectral sequence satisfy the Leibniz rule. Such a result is not known to the author even though it is plausible. If $\alpha$ is a unit, this would imply that the associated formal group is height exactly $n$ in the sense of \cite[Rem. 4.5.7]{Smi11}. In \cite[Rem. 1.3]{LW20}, it is claimed that inverting the image of $v_1$ in $\TC_*(\mathcal{O}_L/\mathbb{S}_{W(\mathbb{F}_q)}[z])$ has the same effect as inverting $\beta$, so this is true at least asymptotically. We speculate that, by choosing a complex orientation carefully, one can arrange that the image of $v_1$ is exactly the element $\beta^d$ considered in loc. cit., where $d=[L(\zeta_p),L]$. 
\end{rem2}
\begin{comment}
\begin{rem2}
Fix $m\ge 2$. Then, since the $K(m)$-localization $L_{K(m)}\K(\mathcal{O}_L)$ of the algebraic K-theory of $\mathcal{O}_L$ is trivial by Mitchell's theorem \cite{Mit90} and  $\TP_*(\mathcal{O}_L/\mathbb{S}_{W(\mathbb{F}_q)}[z])$ is a $\K(\mathcal{O}_L)$-module, we know  that $L_{K(m)}\TP(\mathcal{O}_L/\mathbb{S}_{W(\mathbb{F}_q)}[z])=0$. Since for ring spectra $L_{K(m)}R=0$ if and only if $L_{T(m)}R=0$ by \cite[Lemma 2.3]{LMMT20}, we know that $T(m)\otimes \TC^{-}(\mathcal{O}_L/\mathbb{S}_{W(\mathbb{F}_q)}[z])=T(m)\otimes \TP(\mathcal{O}_L/\mathbb{S}_{W(\mathbb{F}_q)}[z])=0$. Finally, note that $v_m$-multiplication on an $\MU$-module such as $\TP_*(\mathcal{O}_L/\mathbb{S}_{W(\mathbb{F}_q)}[z])$ agrees with $v_m$-multiplication on $F(m)\otimes  E_*^{-}(\mathcal{O}_L/\mathbb{S}_{W(\mathbb{F}_q)}[z])$ on the first tensor factor for any finite spectrum $F(m)$ of type $m$.\footnote{A finite spectrum $F(n)$ has type $n$ if $K(n-1)_*F(n)=0$ and $K(n)_*F(n)\ne 0$.} 
Therefore, we can determine that some power of $v_m$ acts trivially on $\TP_*(\mathcal{O}_L/\mathbb{S}_{W(\mathbb{F}_q)}[z])$, however since $\TP_*(\mathcal{O}_L/\mathbb{S}_{W(\mathbb{F}_q)}[z])$ is polynomial we know that $v_m$-acts trivially and consequently $v_m\in \pi_{2p^m-2}\MU$ maps trivially to $\TC_*^{-}(\mathcal{O}_L/\mathbb{S}_{W(\mathbb{F}_q)}[z])$ and $\TP(\mathcal{O}_L/\mathbb{S}_{W(\mathbb{F}_q)}[z])$. 	
\end{rem2}
\end{comment}

\bibliographystyle{alpha}
\bibliography{ref}

\begin{thebibliography}{BMMS86}

\bibitem[AR08]{AR08}
Christian Ausoni and John Rognes.
\newblock {\em The chromatic red-shift in algebraic K-theory}, volume~40 of
  {\em Monogr. Enseign. Math.}
\newblock L'Enseignement Math\'{e}matique, Geneva, 2008.

\bibitem[BMMS86]{BMMS86}
R.~R. Bruner, J.~P. May, J.~E. McClure, and M.~Steinberger.
\newblock {\em {$H_\infty $} ring spectra and their applications}, volume 1176
  of {\em Lecture Notes in Mathematics}.
\newblock Springer-Verlag, Berlin, 1986.

\bibitem[BMS19]{BMS19}
Bhargav Bhatt, Matthew Morrow, and Peter Scholze.
\newblock Topological {H}ochschild homology and integral {$p$}-adic {H}odge
  theory.
\newblock {\em Publ. Math. Inst. Hautes \'{E}tudes Sci.}, 129:199--310, 2019.

\bibitem[B{\"o}k87]{Bok87a}
M.~B{\"o}kstedt.
\newblock The topological {H}ochschild homology of $\mathbb{Z}$ and of
  $\mathbb{Z}/p\mathbb{Z}$.
\newblock preprint, 1987.

\bibitem[BR05]{BR05}
Robert~R. Bruner and John Rognes.
\newblock Differentials in the homological homotopy fixed point spectral
  sequence.
\newblock {\em Algebr. Geom. Topol.}, 5:653--690 (electronic), 2005.

\bibitem[CM15]{CM15}
Steven~Greg Chadwick and Michael~A. Mandell.
\newblock {$E_n$} genera.
\newblock {\em Geom. Topol.}, 19(6):3193--3232, 2015.

\bibitem[Fr{\"o}68]{Fro68}
A.~Fr{\"o}hlich.
\newblock {\em Formal groups}.
\newblock Lecture Notes in Mathematics, No. 74. Springer-Verlag, Berlin-New
  York, 1968.

\bibitem[HW20]{HW20}
Jeremy {Hahn} and Dylan {Wilson}.
\newblock {Redshift and multiplication for truncated Brown-Peterson spectra}.
\newblock {\em arXiv e-prints}, page arXiv:2012.00864, December 2020.

\bibitem[KN19]{KN19}
Achim {Krause} and Thomas {Nikolaus}.
\newblock {B{\"o}kstedt periodicity and quotients of DVRs}.
\newblock {\em arXiv e-prints}, page arXiv:1907.03477, July 2019.

\bibitem[Laz55]{Laz55}
Michel Lazard.
\newblock Sur les groupes de {L}ie formels \`a un param\`etre.
\newblock {\em Bull. Soc. Math. France}, 83:251--274, 1955.

\bibitem[Lur]{Lur21}
Jacob Lurie.
\newblock {E}lliptic {C}ohomology {II}: {O}rientations.
\newblock Preprint available at
  \url{https://www.math.ias.edu/~lurie/papers/Elliptic-II.pdf}.

\bibitem[LW20]{LW20}
Ruochuan {Liu} and Guozhen {Wang}.
\newblock {Topological Cyclic Homology of Local Fields}.
\newblock {\em arXiv e-prints}, December 2020.

\bibitem[Mil58]{Mil58}
John Milnor.
\newblock The {S}teenrod algebra and its dual.
\newblock {\em Ann. of Math. (2)}, 67:150--171, 1958.

\bibitem[Rav86]{Rav86}
Douglas~C. Ravenel.
\newblock {\em Complex cobordism and stable homotopy groups of spheres}, volume
  121 of {\em Pure and Applied Mathematics}.
\newblock Academic Press, Inc., Orlando, FL, 1986.

\bibitem[Ser79]{Ser79}
Jean-Pierre Serre.
\newblock {\em Local fields}, volume~67 of {\em Graduate Texts in Mathematics}.
\newblock Springer-Verlag, New York, 1979.
\newblock Translated from the French by Marvin Jay Greenberg.

\bibitem[Smi11]{Smi11}
Brian~D. Smithling.
\newblock On the moduli stack of commutative, 1-parameter formal groups.
\newblock {\em J. Pure Appl. Algebra}, 215(4):368--397, 2011.

\end{thebibliography}

\end{document}